\documentclass[english,10pt]{article}
\usepackage[T1]{fontenc}
\usepackage[latin1]{inputenc}
\usepackage{amssymb}

\makeatletter

 \newenvironment{lyxlist}[1]
   {\begin{list}{}
     {\settowidth{\labelwidth}{#1}
      \setlength{\leftmargin}{\labelwidth}
      \addtolength{\leftmargin}{\labelsep}
      }}
   {\end{list}}

\usepackage{babel}
\makeatother
\begin{document}

\begin{center}
{\LARGE \bf The Sum Theorem for Linear Maximal

\smallskip

 Monotone Operators}

\strut

{\Large M.D. Voisei}

{\large Department of Mathematics

The University of Texas -- Pan American, USA

{\small E-MAIL: \texttt{mvoisei@utpa.edu}}}
\end{center}
\strut

\noindent A{\small BSTRACT:} The goal of this article is to give a
positive answer to Rockafellar's maximality of the sum conjecture
in the linear multi-valued operator case.

\strut

\noindent K{\small EY} W{\small ORDS:} \emph{Linear maximal monotone operator, Minkowski
sum}

\strut

\strut

\strut

\noindent \textbf{1. Introduction}. Rockafellar's maximal monotonicity
of the sum conjecture states that

\strut

\textbf{Conjecture.} \emph{If $A,B$ are multi-valued maximal monotone
operators in a Banach space $X$ with\begin{equation}
{\rm int}D(A)\cap D(B)\neq\emptyset,\label{qc}\end{equation}
then $A+B$ is maximal monotone. Here $D(A),D(B)$ stand for the domains
of $A,B$ while {}``int'' denotes the interior.}

\strut

Currently, this problem is still open. The most general additional
assumptions under which the conjecture is known to hold in a non-reflexive
Banach space settings are $A,B$ are subdifferentials or $D(A),D(B)$
are closed convex (see e.g. Simons {[}3{]}, Voisei {[}4{]}).

When $A,B$ are linear the constraint qualification (\ref{qc}) is
equivalent to $D(A)=X$. In this case and for this type of constraint
the conjecture was shown to be true by Phelps \& Simons {[}2{]} for
single-valued operators and later, an improved version, for a different
form of the constraint, was proved in Voisei {[}5{]}.

The aim of the present note is to prove that the conjecture is true
for multi-valued linear maximal monotone operators under a weaker
qualification constraint of the form (\ref{qc}). Our main result
is the following.

\strut

\textbf{Theorem 1.} \emph{Let $A,B:X\rightarrow2^{X^{*}}$ be linear
multi-valued maximal monotone operators in a Banach space $X$. Assume
that $D(A)-D(B)$ is closed in $X$. Then $A+B$ is maximal monotone.}

\strut

For the notations notions and results concerning maximal monotone
operators and convex analysis we refer to the book of Z\u alinescu
{[}6{]} and the references therein.

\strut

\noindent \textbf{2. The Proof of Theorem 1.} Let $(x_{0},x_{0}^{*})$
be monotonically related to $A+B$, that is\begin{equation}
\langle x-x_{0},a^{*}+b^{*}-x_{0}^{*}\rangle\ge0,\label{ineq}\end{equation}
for every $x\in D(A+B):=D(A)\cap D(B)$, $a^{*}\in Ax$, $b^{*}\in Bx$.
Here {}``$\langle\cdot,\cdot\rangle$'' stands for the dual product
between $X$ and its topological dual $X^{*}$.

Consider the convex functions $\alpha,\beta:X\times X^{*}\rightarrow\mathbb{R}\cup\{\infty\}$
and $\Phi:{\cal X}\times{\cal Y}\rightarrow\mathbb{R}\cup\{\infty\}$
given by\[
\alpha(x,x^{*})=\langle x-x_{0},x^{*}\rangle,\ {\rm if}\ (x,x^{*})\in A;\ \alpha(x,x^{*})=+\infty,\ {\rm otherwise},\]
\[
\beta(x,x^{*})=\langle x-x_{0},x^{*}\rangle-\langle x,x_{0}^{*}\rangle,\ {\rm if}\ (x,x^{*})\in B;\ \beta(x,x^{*})=+\infty,\ {\rm otherwise},\]
\[
\Phi(x,x^{*},z^{*};y)=\alpha(x+y,x^{*})+\beta(x,z^{*})+\langle x_{0},x_{0}^{*}\rangle,\]
$(x,x^{*},z^{*})\in{\cal X}:=X\times X^{*}\times X^{*},\  y\in{\cal Y}:=X.$

The convexity of $\alpha,\beta,\Phi$ comes from the fact that $A,B$
are linear monotone. Moreover, $\alpha,\beta,\Phi$ are lower semicontinuous
with respect to the strong topology, since $A,B$ are maximal monotone
thus closed.

Let ${\rm Pr}_{{\cal Y}}(u,v)=v$, $(u,v)\in{\cal X}\times{\cal Y}$,
be the projection of ${\cal X}\times{\cal Y}$ onto ${\cal Y}$. Clearly,
$D(\alpha)=A$, $D(\beta)=B$, and ${\rm Pr}_{{\cal Y}}(D(\Phi))=D(A)-D(B)$.
Because $D(A)-D(B)$ is a closed subspace, it yields that $0\in^{{\rm ic}}({\rm Pr}_{{\cal Y}}D(\Phi))$.
According to {[}6, Theorem 2.7.1 (vii), p. 113{]}, the fundamental
duality formula holds, that is

\begin{equation}
\inf_{u\in{\cal X}}\Phi(u,0)=\max_{y^{*}\in X^{*}}(-\Phi^{*}(0,y^{*})).\label{df}\end{equation}
 Notice that, from (\ref{ineq}), $\Phi(x,x^{*},z^{*};0)=\alpha(x,x^{*})+\beta(x,z^{*})+\langle x_{0},x_{0}^{*}\rangle\ge0$,
for every $u=(x,x^{*},z^{*})\in{\cal X}$, i.e., $\inf\limits _{u\in{\cal X}}\Phi(u,0)\ge0$.
Therefore, (\ref{df}) provides an $y^{*}\in X^{*}$, such that $\Phi^{*}(0;y^{*})=\sup\{\langle y,y^{*}\rangle-\Phi(x,x^{*},z^{*};y);\  x,y\in X,\  x^{*},z^{*}\in X^{*}\}\le0$,
i.e.,\begin{equation}
\alpha(z,x^{*})+\beta(x,z^{*})+\langle x_{0},x_{0}^{*}\rangle-\langle z-x,y^{*}\rangle\ge0,\label{final}\end{equation}
for every $x,z\in X$, $x^{*},z^{*}\in X^{*}$.

We have \[
\inf_{(z,x^{*})\in X\times X^{*}}\alpha(z,x^{*})-\langle z,y^{*}\rangle=\inf_{(z,x^{*})\in A}\langle z-x_{0},x^{*}-y^{*}\rangle-\langle x_{0},y^{*}\rangle,\]
\[
\inf_{(x,z^{*})\in X\times X^{*}}\beta(x,z^{*})+\langle x,y^{*}\rangle=\inf_{(x,z^{*})\in B}\langle x-x_{0},z^{*}-x_{0}^{*}+y^{*}\rangle-\langle x_{0},x_{0}^{*}-y^{*}\rangle.\]
Hence (\ref{final}) is equivalent to\begin{equation}
\inf_{(z,x^{*})\in A}\langle z-x_{0},x^{*}-y^{*}\rangle+\inf_{(x,z^{*})\in B}\langle x-x_{0},z^{*}-x_{0}^{*}+y^{*}\rangle\ge0.\label{infineq}\end{equation}
But, $\inf\limits _{(z,x^{*})\in A}\langle z-x_{0},x^{*}-y^{*}\rangle\le0$
and $\inf\limits _{(x,z^{*})\in B}\langle x-x_{0},z^{*}-x_{0}^{*}+y^{*}\rangle\le0$
because $A,B$ are maximal monotone. Taking (\ref{infineq}) into
account we find\[
\inf_{(z,x^{*})\in A}\langle z-x_{0},x^{*}-y^{*}\rangle=\inf_{(x,z^{*})\in B}\langle x-x_{0},z^{*}-x_{0}^{*}+y^{*}\rangle=0,\]
 that is, $y^{*}\in Ax_{0}$, $x_{0}^{*}-y^{*}\in Bx_{0}$. Therefore
$(x_{0},x_{0}^{*})\in A+B$. The proof is complete. $\square$

\strut

\textbf{Remark 1.} As seen in {[}2, Remark 7.3, p. 326{]} or in {[}3,
Problem 34.2, p. 104{]} the condition $D(A)-D(B)$ closed cannot be
further relaxed.

\strut

\textbf{Remark 2.} An alternative proof for Therem 1 can be given by reducing the condition $D(A)-D(B)$ closed to $D(A)-D(B)=X$, based on the fact that the maximality of any linear monotone operator $S$ in $X\times X^*$ is equivalent to the maximality of $S$ in $Y\times Y^*$, for every subspace $Y$ of $X$ that contains $D(S)$. In this case the $Y-$saturation condition is easily satisfied (see e.g. [3, Theorem 16.10, p. 44]).

\strut

\textbf{Acknowledgements.} The author would like to thank Dr. G. Moro\c{s}anu for interesting suggestions and comments with regard to a preliminary version of this paper.

\strut

\strut

\centerline{ \textbf{References}}

\begin{lyxlist}{00.00.0000}
\item [{[}1{]}]S. Fitzpatrick, \emph{Representing monotone operators by
convex functions}, Workshop/Miniconference on Functional Analysis
and Optimization (Canberra, 1988), 59--65, Proc. Centre Math. Anal.
Austral. Nat. Univ., 20, Austral. Nat. Univ., Canberra, 1988. 
\item [{[}2{]}]R.R. Phelps, S. Simons, \emph{Unbounded linear monotone
operators on nonreflexive Banach spaces (English. English summary)},
J. Convex Anal. \textbf{5} (2) (1998), 303--328.
\item [{[}3{]}]S. Simons, \emph{Minimax and monotonicity}, Lecture Notes
in Mathematics, 1693. Springer-Verlag, Berlin, 1998.
\item [{[}4{]}]M.D. Voisei, \emph{A maximality theorem for the sum of maximal
monotone operators in non-reflexive Banach spaces}, Math. Sci. Res.
J., \textbf{10} (2) 2006, 36--41.
\item [{[}5{]}]M.D. Voisei, \emph{Monotonicity representability and maximality
via the Fitzpatrick function}, preprint 2006.
\item [{[}6{]}]C. Z\u alinescu, \emph{Convex analysis in general vector
spaces}, World Scientific Publishing Co., Inc., River Edge, NJ, 2002.
\end{lyxlist}

\end{document}